\documentclass[11pt]{article}
\usepackage{amsmath, amscd, amssymb, latexsym, epsfig, color, amsthm}
\numberwithin{equation}{section}

\newtheorem{theorem}{Theorem}
\newtheorem{proposition}[theorem]{Proposition}

\newtheorem{lemma}[theorem]{Lemma}

\theoremstyle{definition}

\DeclareMathOperator{\conv}{\mathrm{conv}}

\newcommand{\R}{{\mathbb R}}

\newcommand{\fmax}{\ensuremath{\mathrm{fmax}}\hspace{1pt}}

\newcommand{\Int}{\mbox{\upshape int}\,}

\title{From acute sets to centrally symmetric $2$-neighborly polytopes}
\author{Isabella Novik\thanks{Research is partially\textsl{}
supported by NSF grant DMS-1664865 and by Robert R.~\&  Elaine F.~Phelps Professorship in Mathematics. 
}\\
\small Department of Mathematics\\[-0.8ex]
\small University of Washington\\[-0.8ex]
\small Seattle, WA 98195-4350, USA\\[-0.8ex]
\small \texttt{novik@math.washington.edu}
}

\begin{document}
\maketitle

\begin{abstract} What is the maximum number of vertices that a centrally symmetric 2-neighborly polytope of dimension $d$ can have? It is known that the answer does not exceed $2^d$. Here we provide an explicit construction showing that it is at least $2^{d-1}+2$.
\end{abstract}

\section{Introduction}
The goal of this note is to construct centrally symmetric $2$-neighborly polytopes with many vertices. Recall that a polytope is the convex hull of a set of finitely many points in $\R^d$. The dimension of a polytope $P$ is the dimension of its affine hull. We say that $P$ is a $d$-polytope if the dimension of $P$ is equal to $d$. A polytope $P$ is {\em centrally symmetric} (cs, for short) if for every $x\in P$, $-x$ belongs to $P$ as well. 

A cs polytope $P$ is called {\em $k$-neighborly} if every set of $k$ of its vertices, no two of which are antipodes, is the vertex set of a face of $P$. In addition to being of intrinsic interest, the study of cs $k$-neighborly polytopes is motivated by the recently discovered tantalizing connections (initiated by Donoho and his collaborators \cite{Don06a,DonTan10}) between such polytopes and seemingly distant areas of error-correcting codes and sparse signal reconstruction. It is also worth mentioning that in contrast with the situation for polytopes without a symmetry assumption, a cs $d$-polytope with sufficiently many vertices cannot be even $2$-neighborly \cite{Burton, LinNov}. 

A few more definitions are in order. A set $S\subset \R^d$ is {\em acute} if every three points from $S$ determine an acute triangle. A set $S\subset \R^d$ is {\em antipodal} if for every two points $x,y\in S$, there exist two (distinct) parallel hyperplanes $H_x$ and $H_y$ such that $x\in H_x$, $y\in H_y$, and all elements of $S$ lie in the closed strip defined by $H_x$ and $H_y$. 

It is well-known and easy to check that any acute set is antipodal. Similarly, as was observed in \cite[Lemma 2.1]{LinNov}, the vertex set of any cs $2$-neighborly polytope $P\subset \R^d$ is antipodal. Furthermore, according to a celebrated theorem of Danzer and Gr\"unbaum  \cite{DanzGrunb} (see also \cite[Ch.~17]{AigZieg}), an antipodal subset of $\R^d$ has at most $2^d$ elements, and it has exactly $2^d$ elements if and only if it is the vertex set of a parallelepiped. Since the vertex set of a parallelepiped is not acute and since a $d$-parallelepiped (for $d>2$) is not cs $2$-neighborly, it follows that any acute set $S\subset\R^d$ has at most $2^d-1$ elements, while any cs $2$-neighborly $d$-polytope $P$ with $d\geq 3$ has at most $2^d-2$ vertices.

Although the size of the largest acute set in $\R^d$ remains a mystery, in a very recent breakthrough paper \cite{GerHar}, Gerencs\'er and Harangi constructed an acute set in $\R^d$ of size $2^{d-1}+1$. The previous record size was $F_{d+2}=\Theta\big((\frac{1+\sqrt{5}}{2})^d\big)$, where $F_n$ is the $n$-th Fibonacci number, see \cite{Zakh}. 

Similarly, the current record size of the vertex set of a cs $2$-neighborlly $d$-polytope is about $\sqrt{3}^d$, \cite{BarvLeeN-antipodal}. We modify the construction of Gerencs\'er and Harangi to establish the following result.

\begin{theorem} \label{main-thm}
There exists a cs $2$-neighborly $d$-polytope with $2^{d-1}+2$ vertices.
\end{theorem}

We say that a (finite) set $S\subset \R^d\setminus\{0\}$ is {\em cs} if for every $x\in S$, the point $-x$ is also in $S$. Observe that a cs set can never be acute: indeed, for any $x,y\in\R^d$, the parallelogram determined by $x, y, -x, -y$ has a non-acute angle. The main insight of this note is the notion of an almost acute set: a set $S$ is {\em almost acute} 
if it is cs and for every ordered triple $(x,y,z)$ of distinct points in $S$, the angle $\angle xyz$ is acute as long as $x$ and $z$ are not antipodes. 

With this definition in hand, the following two results yield Theorem \ref{main-thm}.
\begin{lemma} \label{cs-2-neighb}
Let $S\subset \R^d$ be a cs set that spans $\R^d$, and let $P=\conv(S)$. If $S$ is almost acute, then $P$ is a cs $2$-neighborly $d$-polytope whose vertex set is $S$.
\end{lemma}
\begin{lemma} \label{large-almost-acute}
There exists an almost acute subset of $\R^d$ of size $2^{d-1}+2$.
\end{lemma}
To prove Lemma \ref{large-almost-acute} we modify the Gerencs\'er--Harangi construction: as in \cite{GerHar}, we start with the vertex set of the $(d-1)$-cube $[-1,1]^{d-1}$ embedded in the coordinate hyperplane $\R^{d-1}\times\{0\}$ of $\R^d$. We then use the extra dimension to perturb the vertices in such a way that the resulting set in $\R^d$ is almost acute. (In particular, any pair of antipodes is perturbed to a pair of antipodes.) Adding to this set a pair of antipodes of the form $(0,\ldots,0, c)$ and $(0,\ldots,0, -c)$, where $c\in\R$ is sufficiently large, completes the construction.

The proofs of Lemmas \ref{cs-2-neighb} and \ref{large-almost-acute} are given in Sections 2 and 3, respectively. We close in Section 4 with some remarks and open problems.

\section{Polytopes with an almost acute vertex set}
The goal of this section is to prove Lemma \ref{cs-2-neighb}. For all undefined terminology pertaining to polytopes, we refer our readers to Ziegler's book \cite{Ziegler}. Thus, assume that $S\subset\R^d$ is an {\em almost acute} set that spans $\R^d$. Then $P:=\conv(S)$ is a cs $d$-polytope whose vertex set is contained in $S$. To prove the lemma, we have to show that (i) every $x\in S$ is a vertex of $P$, and (ii) for every $x,y\in S$ with $y\notin\{x,-x\}$, the line segment $[x,y]$ is an edge of $P$.

Let $x$ be any element of $S$. Let $H$ be the hyperplane that contains $x$ and is perpendicular to the line segment $[-x,x]$. Since $S$ is an almost acute set, for every $y\in S\setminus\{-x,x\}$, the angle $\angle (-x)xy$ is acute, so that $y$ lies in the same open half-space of $\R^d$ defined by $H$ as $-x$. It follows that $H$ is a supporting hyperplane of $P$ and that $H\cap P=\{x\}$. Hence $x$ is a vertex of $P$.

Now, let $x,y$ be any elements of $S$ with $y\notin\{x,-x\}$. Consider parallelogram $Q$ with vertices $x,y,-x,-y$. There are two possible cases:

\medskip\noindent{\bf Case 1:} $Q$ is a {\bf rectangle}. Let $H$ be the hyperplane perpendicular to the line segment $[-y,x]$ and passing through $x$, and hence also through $y$. Since $S$ is an almost acute set, for every $z\in S\setminus\{x, -y, y\}$, the angle $\angle(-y)xz$ is acute, so that $z$ lies in the same open half-space of $\R^d$ defined by $H$ as $-y$. We conclude that $H$ is a supporting hyperplane of $P$ and that $H\cap P=\conv(x,y)=[x,y]$. Thus $[x,y]$ is an edge of $P$. 

\medskip\noindent{\bf Case 2:} $Q$ is {\bf not a rectangle}. In this case exactly one of the angles $\angle (-y)xy$, $\angle (-x)yx$ is obtuse. Without loss of generality (by switching the roles of $x$ and $y$ if necessary), we may assume that $\angle (-y)xy$ is obtuse. As in Case 1, let $H$ be the hyperplane perpendicular to the line segment $[-y,x]$ and passing through $x$. Our assumption that $S$ is almost acute then yields that all elements of $S\setminus \{x,y, -y\}$ are contained in the same (open) side of $H$ as $-y$, while $y$ lies on the opposite side of $H$. Therefore, either $[x,y]$ is an edge of $P$, in which case we are done, or all neighbors of $x$ in $P$ lie on the side of $H$ that does not contain $y$. In this latter case, the cone based at $x$ and spanned by the rays from $x$ to the neighbors of $x$ does not contain $y$. This however contradicts the well-known fact (see \cite[Lemma 3.6]{Ziegler}) that such a cone must contain $P$, and hence also $y$. The lemma follows.

\section{Construction of an almost acute set} 
In this section we prove Lemma \ref{large-almost-acute}. To do so, we constract an almost acute set in $\R^d$ of size $2^{d-1}+2$. We start with the set $S^0$ described below; we then perturb the points of $S^0$ to obtain an almost acute set. As our construction/proof is a simple modification of that in \cite{GerHar}, we only sketch the main ideas leaving out some of the details.

Pick a real number $c>\sqrt{d-1}$ and consider the following subset of $\R^d$ of size $2^{d-1}+2$:
$$S^0:=\big\{(\delta_1,\ldots,\delta_{d-1},0) \mid \delta_1,\ldots,\delta_{d-1}\in \{\pm 1\} \big\} \cup \big\{(0,\ldots,0,\pm c)\big\}.$$
Thus, $S^0$ consists of the vertex set $V^{0}$ of the $(d-1)$-cube $[-1,1]^{d-1}\times\{0\}\subset \R^{d}$
 and two additional points, $x_0$ and $-x_0$, positioned high above and far below the center of the cube, respectively. An easy computation shows that for all distinct $y,z\in V^{0}$, the angles $\angle(\pm x_0)yz$, $\angle y(\pm x_0)z$, $\angle(\pm x_0)(\mp x_0)y$ are acute, and, assuming also that $z\neq -y$, so is $\angle y(-y)z$. Hence there exists an $\epsilon_0>0$ such that if all vertices of the cube are perturbed by no more than $\epsilon_0$, then all of the above angles remain acute.  Therefore, to complete the proof, it suffices to perturb the vertices of $V^{0}$ in such a way that (i) antipodes are perturbed to antipodes, and (ii) for all $x,y,z\in V^{0}$ no two of which are antipodes, the perturbed triangle is acute.

The key fact we will use is the following lemma from \cite{GerHar}:
\begin{lemma} \label{aux}
Let $V^{0}$ be the vertex set of the $(d-1)$-cube  $[-1,1]^{d-1}\times\{0\}\subset \R^{d}$. For every $\epsilon>0$ and $x\in V^{0}$, there exists $x'\in\R^d$ such that $x'$ is within distance $\epsilon$ from $x$ and the angles $\angle x'yz$ and $\angle yx'z$ are acute for all $y,z \in  V^0\setminus\{x\}$.
\end{lemma}

Arbitrarily order the elements of $V^{0}$, so that  $V^0=\{x_1,-x_1,x_2,-x_2,\ldots, x_{2^{d-2}}, -x_{2^{d-2}}\}$. We induct on $1\leq p\leq 2^{d-2}$, to construct a set  $V^p=\{x'_1,-x'_1, \ldots, x'_p,-x'_p\} \cup \{x_j, -x_j \mid p<j\leq 2^{d-2}\}$ with the property that (a) for all $1\leq i\leq p$, $\|x'_i-x_i\|<\epsilon_0$, and (b) for every three points $x, y, z$ of $V^{p}$ no two of which are antipodes and such that $x=\pm x'_i$ for some $1\leq i\leq p$, the angles $\angle xyz$ and $\angle yxz$ are acute. We refer to (a) and (b) combined as the $(*_p)$-property.

Assume $V^{p-1}$ satisfies the $(*_{p-1})$-property. In particular, for every three points $x'_i, y, z$ of $V^{p-1}$ no two of which are antipodes, the angles $\angle (\pm x'_i)yz$ and $\angle y(\pm x'_i)z$ are acute. Hence there exists an $0<\epsilon_p<\epsilon_0$, such that if $x_p$ and $-x_p$ are perturbed by no more than $\epsilon_p$, then all of the above angles involving $\pm x_p$ (as $y$ or $z$) remain acute. Furthermore, by Lemma \ref{aux}, there exists $x'_p$ within distance $\epsilon_p$ of $x_p$ such that for all $y,z\in\{x_j, -x_j \mid p<j\leq 2^{d-2}\}$, the angles $\angle x'_pyz$ and $\angle yx'_pz$ are acute. Since the set $\{x_j, -x_j \mid p<j\leq 2^{d-2}\}$ is cs, it follows that the angles $\angle(-x'_p)yz$ and $\angle y(-x'_p)z$ are also acute. We conclude that $V^p$ satisfies the $(*_p)$-property. The set $S:=V^{2^{d-2}}\cup \{x_0,-x_0\}$ is then an almost acute set of size $2^{d-1}+2$. This completes the proof of Lemma \ref{large-almost-acute}, and hence also of Theorem \ref{main-thm}.

\section{Concluding remarks and open problems}
We close with a few open problems.

The main result of this note together with \cite[Lemma 2.1]{LinNov} and the Danzer--Gr\"unbaum theorem \cite{DanzGrunb} implies that for $d\geq 3$, the maximum number of vertices that a cs $2$-neighborly $d$-polytope can have lies in the interval $[2^{d-1}+2, 2^d-2]$. In dimension three, the only cs $2$-neighborly polytope is the cross-polytope, which indeed has $6=2^2+2=2^3-2$ vertices. In dimension four, the maximum is $10=2^3+2$; this result is due to Gr\"unbaum, see \cite[p.~116]{Gru-book}. However, for $d>4$, the exact value of the maximum remains unknown. 

A related question is what is the maximum number of edges, $\fmax(d;N)$, that a cs $d$-polytope with $N$ vertices can have. At present, it is known that 
\begin{equation} \label{edges}
\left(1-3^{-\lfloor d/2-1\rfloor}\right){N \choose 2}\leq \fmax(d;N) \leq \left(1-2^{-d}\right)\frac{N^2}{2};
\end{equation}
see \cite[Theorem 3.2(2)]{BarvLeeN-antipodal} and \cite[Proposition 2.1]{BarvN} for the lower and the upper bound, respectively. However, the exact value of $\fmax(d;N)$ or even its asymptotics remains a mystery. The main result of this paper makes us believe that $\fmax(d;N)$ might be closer to the right-hand side of Eq.~\eqref{edges} than to the left one.

Finally, it would be interesting to understand the maximum number of vertices that a cs $3$-neighborly $d$-polytope can have. It is known that there exist cs $3$-neighborly $d$-polytopes with $\approx 2^{0.023d}$ vertices, see \cite[Remark 4.3]{BarvLeeN-antipodal}. On the other hand, an argument similar to the proof of \cite[Theorem 1.1]{LinNov}, shows that a cs $d$-polytope with $\lceil 2\sqrt{2}\cdot 3^{0.5d}\rceil$ or more vertices cannot be $3$-neighborly.

\section*{Acknowledgements} The author is grateful to a participant of the mathematical forum dxdy.ru, who wished to remain anonymous, for bringing reference \cite{GerHar} to our attention.
{\small
\bibliography{refs}
\bibliographystyle{plain}
}
\end{document}